\documentclass[10pt, twoside]{article}
\usepackage{amssymb,latexsym,amsmath}
\usepackage[hypertex]{hyperref}
\usepackage[english]{babel}
\begin{document}
\makeatletter
\renewcommand{\@evenfoot}{ \thepage \hfil \footnotesize{\it ISSN 1025-6415 \ \ Reports of the National Akademy of Sciences of Ukraine, 
2003, no. 8} } 
\renewcommand{\@oddfoot}{\footnotesize{\it ISSN 1025-6415 
\ \  Dopovidi Natsionalno\" \i \ Akademi\" \i \  Nauk Ukra\"\i ni, 
2003, no. 8} \hfil \thepage } 
\noindent
\\
\\
\\
\\
\\
\\
\\
\\
\\
\par{
\leftskip=1.5cm  \rightskip=0cm  
\noindent UDC 512\medskip \\
\copyright \ {\bf 2003}\medskip \\
{\bf T. R. Seifullin}\bigskip \\ 
{\large \bf
Determination of the basis of the space of all root \medskip\\
functionals of a system of polynomial equations \medskip\\
and of the basis of its ideal by the operation of the extension\medskip\\ 
of bounded root functionals \bigskip \\
} 
{\it 
(Presented by Corresponding Member of the NAS of Ukraine A. A. Letichevsky)
}\medskip \\
\small {\it 
The notion of a root functional of a system of polynomials or
an ideal of polynomials is
a generalization of the notion of a root, in particular, for a multiple root.
A basis of the space of all root functionals
and a basis of the ideal are found
by using the operation of extension of bounded root functionals
when the number of equations is equal to the number of unknowns
and if it is known that the number of roots is finite.
The asyptotic complexity of these methods is~$d^{O(n)}$ operations,
where $n$ is the number of equations and unknowns, $d$ is the maximal degree
of polynomials.}
\par \bigskip } \noindent 
Presence of roots at infinity leads to large degrees 
 of polynomials  in Buchberger algorithm 
for construction of a Gr\"obner basis of the ideal of polynomials [8].  
Therefore the complexity of Buchberger algorithm 
such large, in the case of the $0$-dimensional  variety of roots it 
is equal to $d^{O(n^2)}$ for the number of operations [9], where $d$ 
is the maximal degree of polynomials, $n$~is the number of variables. 
In the paper [10] it is shown the exactness of this estimation. 
\hbox{For a system} of polynomial equations, in which the number of polynomials 
is equal to the number of variables, the application 
of extension operations to bounded  root functionals [6], [7] 
gradually cuts components of \hbox{functionals,} lying at infinity, 
not exiting over the limits of degrees 
$\leq (d_1-1)+\ldots +(d_n-1)$, where $d_1,\ldots ,d_n$  
are degrees of polynomials. 
This allows, in the case, if it is known, that the \hbox{variety} 
of roots is $0$-dimensional, to find a basis of the space 
of all root functionals of the system of polynomials and 
a basis of the ideal of polynomials in $O(D^4)$  
operations, where $D=C^n_{d_1+\ldots +d_n}$. 
A similar complexity is 
had 
by the method, based  on the use of a 
multivariate resultant, 
that find all \hbox{isolated} roots  of polynomials  in $d^{O(n)}$ operations, 
even in the case of the infinite number of roots at affine domain 
and at infinity [11]. \medskip

 Let ${\bf R} $ a be commutative ring with unity $1$ and zero 
$0$.

 Let $x=(x_1,\ldots ,x_n)$  be variables, 
${\bf R} [x]$ be a ring of polynomials in variables $x$ 
with coefficients in ${\bf R} $.

 In the paper we will use definition and assumption, given 
in [6,7]. \smallskip  {\it

 {\bf Lemma 1.}   Let  $x=(x_1,\ldots ,x_n)$   be  variables, 
 $f(x)=(f_1(x),\ldots ,f_s(x))$   be polynomials. There holds:

 1) a functional $L(x_*)$ annuls 
$(f(x))_x$ if and only if $\forall i=1,s: L(x_*)\cdot f_i(x) = 0$;

 2) a functional $L(x_*)$ annuls \ $(f(x))^{\leq d} _x$ \   
if and only if \ $\forall i=1,s: L(x_*)\cdot f_i(x) = 0$ in 
${\bf R} [x^{\leq d-\deg (f_i)} ]$. } \smallskip

 {\bf Proof 1.} \ $L(x_*)$ \ annuls \ $(f(x))_x$ \ if  and  only if \ 
$\forall i=1,s:  0  =  L(x_*).f_i(x)\cdot g^i(x) = $
\hbox{$L(x_*)\cdot f_i(x).g^i(x)$} 
for any $g^i(x) \in {\bf R} [x]$. 
$L(x_*)\cdot f_i(x).g^i(x) = 0$ for any $g^i(x) \in {\bf R} [x]$ 
means that $L(x_*)\cdot f_i(x) = 0$.

 {\bf Proof 2.}  $L(x_*)$ annuls $(f(x))^{\leq d} _x$ 
if and only if $\forall i=1,s: 0=L(x_*).f_i(x)\cdot g^i(x)=$
\hbox{$L(x_*)\cdot f_i(x).g^i(x)$} 
for  any  $g^i(x) \in {\bf R} [x^{\leq d-\deg (f_i)} ]$. 
$L(x_*)\cdot f_i(x).g^i(x) = 0$ 
for  any  $g^i(x)  \in   {\bf R} [x^{\leq d-\deg (f_i)} ]$ 
means that $L(x_*)\cdot f_i(x) = 0$ in 
${\bf R} [x^{\leq d-\deg (f_i)} ]$. \smallskip

 {\bf Definition 1.}  Let ${\cal V} $  be a module  over  
${\bf R} $,  denote by  ${\cal V} _*$  the set of all linear 
over ${\bf R} $ maps ${\cal V}  \rightarrow  {\bf R} $. 
Let ${\cal U} $  be a submodule of the module ${\cal V} $ over ${\bf R} $, 
denote by ${\cal U} ^\bot $ the set of all $l\in {\cal V}_*$, 
annulling ${\cal U} $, i. e. such  that $\forall F\in 
{\cal U} : l.F=0$. \smallskip

 {\bf Definition 2.}  Let ${\cal U} ,{\cal V} ,{\cal G} 
$  be sets, let  $l:  {\cal V}   \rightarrow   {\cal G} $  
be a map, let ${\cal U} \subseteq {\cal V} $. 
Denote by $l| _{\cal U} $ the restriction of the map $l$ on 
the set ${\cal U} $, i. e. such a map $l': {\cal U}  \rightarrow  {\cal 
G} $,  that $\forall F\in {\cal U} : l'.F = l.F$. \smallskip

{\it {\bf Statement 1.}  Let ${\cal U} $  be a submodule of a module 
${\cal V} $ over ${\bf R} $, let  ${\cal L} $   be a submodule 
of the module ${\cal V} _*$ over ${\bf R} $. If $l_1,l_2 \in  {\cal 
L} $, then $l_1 = l_2$ in ${\cal U} $ if and only if 
$l_1 - l_2 \in  {\cal U} ^\bot $, the last means that  
$l_1/{\cal U} ^\bot   =  l_2/{\cal U} ^\bot $.  Hence,  there is  
an isomorphism ${\cal L} | _{\cal U}  \simeq  {\cal L} /{\cal U} 
^\bot $ such that 
$l| _{\cal U} \leftrightarrow l/{\cal U} ^\bot $ 
for  any  $l  \in   {\cal L} $. \smallskip } {\it

 {\bf Theorem 1.}  Let  $x=(x_1,\ldots ,x_n)$   be  variables, 
 $f(x)=(f_1(x),\ldots ,f_n(x))$   be polynomials, ${\delta} _f = 
\sum\limits^{ n} _{i=1} (\deg (f_i)-1)$.

 1. Let ${\delta} ,{\delta} '\geq 0$. If  $L(x_*)$  annuls 
 $(f(x))^{\leq {\delta} _f+{\delta} } _x$,  $L'(x_*)=L(x_*)$ 
 in ${\bf R} [x^{\leq {\delta} _f+{\delta} '} ]$, then $L'(x_*)$ 
annuls $(f(x))^{\leq {\delta} _f+{\delta} } _x\cap {\bf 
R} [x^{\leq {\delta} _f+{\delta} '} ]$.

 2.  Let  $0\leq {\delta} '_1\leq {\delta} _1$,  $0\leq {\delta} 
'_2\leq {\delta} _2$,  $0\leq {\delta} \leq {\delta} _1+{\delta} 
_2+1$,  $0\leq {\delta} '\leq {\delta} '_1+{\delta} '_2+1$.  
Let  $L_1(x_*)$ annuls $(f(x))^{\leq {\delta} _f+{\delta} 
_1} _x$, \ $L_2(x_*)$ annuls $(f(x))^{\leq {\delta} _f+{\delta} 
_2} _x$, \ then $L_1(x_*)*L_2(x_*)$ annuls $(f(x))^{\leq 
{\delta} _f+{\delta} } _x$. If $L'_1(x_*) = L_1(x_*)$ in ${\bf 
R} [x^{\leq {\delta} _f+{\delta} '_1} ]$, $L'_2(x_*)  =  L_2(x_*)$ 
in ${\bf R} [x^{\leq {\delta} _f+{\delta} '_2} ]$, then $L'_1(x_*)*L'_2(x_*) 
= L_1(x_*)*L_2(x_*)$ in ${\bf R} [x^{\leq {\delta} _f+{\delta} 
'} ]$.

 From above it follows that the extension map $*$ for functionals
induces the map \medskip

\ \ \ $\frac{ ((f(x))^{\leq 
{\delta} _f+{\delta} _1} _x)^\bot } {  {\bf R} [x^{\leq {\delta} 
_f+{\delta} '_1} ]^\bot }  \times  \frac{ ((f(x))^{\leq {\delta} 
_f+{\delta} _2} _x)^\bot } {  {\bf R} [x^{\leq {\delta} _f+{\delta} 
'_2} ]^\bot }  \rightarrow  \frac{ ((f(x))^{\leq {\delta} _f+{\delta} 
} _x)^\bot } {   {\bf R} [x^{\leq {\delta} _f+{\delta} '} ]^\bot },$ 
\medskip \\  
or,  in other words, induces the map \medskip

\ \ \ $((f(x))^{\leq 
{\delta} _f+{\delta} _1} _x)^\bot | _{{\bf R} [x^{\leq {\delta} 
_f+{\delta} '_1} ]}     \times      ((f(x))^{\leq {\delta} _f+{\delta} 
_2} _x)^\bot | _{{\bf R} [x^{\leq {\delta} _f+{\delta} '_2} ]} 
     \rightarrow ((f(x))^{\leq {\delta} _f+{\delta} } _x)^\bot 
| _{{\bf R} [x^{\leq {\delta} _f+{\delta} '} ]}. $ } \eject

 {\bf Proof  1.} Let 
$F(x)\in (f(x))^{\leq {\delta} _f+{\delta} } _x\cap 
{\bf R} [x^{\leq {\delta} _f+{\delta} '} ]$,  
$L'(x_*).F(x) =L(x_*).F(x)$, since $L'(x_*)=L(x_*)$ 
in ${\bf R} [x^{\leq {\delta} _f+{\delta} '} ] \ni  F(x)$ 
; and 
$L(x_*).F(x)=  0$,  since  $L(x_*)$  annuls 
 $(f(x))^{\leq {\delta} _f+{\delta} } _x  \ni   F(x)$.  Then, 
 by the arbitrariness of   
$F(x) \in (f(x))^{\leq {\delta} _f+{\delta} } _x\cap 
{\bf R} [x^{\leq {\delta} _f+{\delta} '} ]$,    
$L'(x_*)$  annuls 
\hbox{$(f(x))^{\leq {\delta} _f+{\delta} } _x
\cap {\bf R} [x^{\leq {\delta} _f+{\delta} '} ]$.}

 {\bf Proof 2.}  Since \ $L_1(x_*)$ \ annuls 
 $(f(x))^{\leq {\delta} _f+{\delta} _1} _x$,  $L_2(x_*)$ annuls 
$(f(x))^{\leq {\delta} _f+{\delta} _2} _x$, then by virtue of {\it 2} 
of theorem  3  in  [6]  $L_1(x_*)*L_2(x_*)$ annuls $(f(x))^{\leq 
{\delta} _f+{\delta} _1+{\delta} _2+1} _x  \supseteq  (f(x))^{\leq 
{\delta} _f+{\delta} } _x$,  hence,  annuls $(f(x))^{\leq 
{\delta} _f+{\delta} } _x$.

 Since $L_1(x_*)$ annuls 
$(f(x))^{\leq {\delta} _f+{\delta} _1} _x \supseteq  
(f(x))^{\leq {\delta} _f+{\delta} '_1} _x$, 
then  annuls and 
$(f(x))^{\leq {\delta} _f+{\delta} '_1} _x$; 
since $L_2(x_*)$ annuls $(f(x))^{\leq {\delta} _f+{\delta}_2} _x 
\supseteq  (f(x))^{\leq {\delta} _f+{\delta} '_2} _x$, 
then annuls and $(f(x))^{\leq {\delta} _f+{\delta} '_2} _x$. 
Then, since $L'_1(x_*) = L_1(x_*)$ in 
${\bf R} [x^{\leq {\delta} _f+{\delta} '_1} ]$  
and $L'_2(x_*) = L_2(x_*)$ in ${\bf R} [x^{\leq {\delta} _f+{\delta} '_2} ]$, 
by virtue of {\it 3} of theorem 3 in  [6]  
\hbox{$L'_1(x_*)*L'_2(x_*)=$} $L_1(x_*)*L_2(x_*)$ in 
${\bf R} [x^{\leq {\delta} _f+{\delta} '_1+{\delta} '_2+1} ] \supseteq 
 {\bf R} [x^{\leq {\delta} _f+{\delta} '} ]$, hence, 
 $L'_1(x_*)*L'_2(x_*)= L_1(x_*)*L_2(x_*)$ in ${\bf R} [x^{\leq 
{\delta} _f+{\delta} '} ]$.  Two last  statements  are obtained 
 by applying of statement 1. \smallskip

 {\bf Definition 3.}  Let $x=(x_1,\ldots ,x_n)$  be variables, \ 
$f(x)=(f_1(x),\ldots ,f_n(x))$  be polynomials, $d\geq 0$. Denote by 
${\cal P} ^{\leq d} _x$ a linear over ${\bf R} $ map  
${\bf R} [x]  \rightarrow   {\bf R} [x]$ such that ${\cal 
P} ^{\leq d} _x.x^{\alpha}  = x^{\alpha} $, if $| {\alpha} 
| \leq d$, and ${\cal P} ^{\leq d} _x.x^{\alpha}  = 0$, if $| 
{\alpha} | >d$. {\it

 {\bf Statement 2.}  Let $x=(x_1,\ldots ,x_n)$  be variables, 
$f(x)=(f_1(x),\ldots ,f_n(x))$  be polynomials, let $d\geq 0$.

 If $L(x_*) \in  {\bf R} [x]_*$, then $L(x_*) = L(x_*).{\cal 
P} ^{\leq d} _x$ in ${\bf R} [x^{\leq d} ]$ and  $L(x_*).{\cal 
P} ^{\leq d} _x  =  0$  in ${\bf R} [x^{>d} ]$.

 If $l(x_*) \in  {\bf R} [x^{\leq d} ]_*$, then the functional 
$L(x_*) = l(x_*).{\cal P} ^{\leq d} _x \in  {\bf R} [x]_*$ and 
 is continuation of $l(x_*)$ on ${\bf R} [x]$, i.  e.  $l(x_*) 
 =  L(x_*)$  in ${\bf R} [x^{\leq d} ]$,  moreover,  $L(x_*).{\cal 
P} ^{\leq d} _x = 0$ in ${\bf R} [x^{>d} ]$. }

 {\bf Proof.}  Let $| {\alpha} | \leq d$, then 
 $L(x_*).{\cal P} ^{\leq d} _x.x^{\alpha}   =  L(x_*).x^{\alpha} 
$.  Since the monoms $x^{\alpha} $, for which $| {\alpha} | \leq 
d$,  linearly  over  ${\bf R} $  generate  ${\bf R} [x^{\leq d} 
]$,  then  $L(x_*)  =L(x_*).{\cal P} ^{\leq d} _x$ in ${\bf 
R} [x^{\leq d} ]$.

 Let $| {\alpha} | >d$, then \ $L(x_*).{\cal P} ^{\leq d}_x.x^{\alpha} = 0$.  
Since the monoms  $x^{\alpha} $, for which $| {\alpha} | >d$, 
linearly over ${\bf R} $ generate ${\bf R} [x^{>d} ]$, 
then $L(x_*).{\cal P} ^{\leq d} _x = 0$ in ${\bf R} [x^{>d} ]$.

 The second part of the statement proved exactly as the first. 

 {\bf Commentary to theorem 1.}  In {\it 2} of theorem 1  computation of 
 $L'_1(x_*)*L'_2(x_*)$ on ${\bf R} [x^{\leq {\delta} _f+{\delta} 
'} ]$  use values of $L'_1(x_*)$ outside ${\bf R} [x^{\leq 
{\delta} _f+{\delta} '_1} ]$  and  values of  $L'_2(x_*)$ outside 
${\bf R} [x^{\leq {\delta} _f+{\delta} '_2} ]$, therefore necessary 
to determine values of $L'_1(x_*)$  outside  ${\bf R} [x^{\leq {\delta} 
_f+{\delta} '_1} ]$ and values of $L'_2(x_*)$ 
\hbox{outside} ${\bf R} [x^{\leq 
{\delta} _f+{\delta} '_2} ]$. With computational point of view 
it is convenient to determine the action of 
$L'_1(x_*)$ in ${\bf R} [x^{>{\delta} _f+{\delta} 
'_1} ]$, and the action of $L'_2(x_*)$ in ${\bf R} [x^{>{\delta} _f+{\delta} 
'_2} ]$  as  zeroes. This holds in the case, if we set 
 $L'_1(x_*)  =  L_1(x_*).{\cal P} ^{\leq {\delta} _f+{\delta} 
'_1} _x$,  $L'_2(x_*)  =L_2(x_*).{\cal P} ^{\leq {\delta} 
_f+{\delta} '_2} _x$. It is enough to compute values 
of the functional $L'_1(x_*)*L'_2(x_*)$ only on 
${\bf R} [x^{\leq {\delta} _f+{\delta} '} ]$. 

 {\bf Definition  4.}   Let  $x=(x_1,\ldots ,x_n)$,   $y\simeq 
x$    be   variables,   $f(x)   =(f_1(x),\ldots ,f_n(x))$ 
 be polynomials.  A functional  $E(x_*)$  we call a unit  root 
functional  of polynomials  $f(x)$,  if  it  annuls  $(f(x))_x$, 
 and  $E(x_*)*1   =E(y_*).\det \left\| \nabla f(x,y)\right\| 
 =  1  +  f(x)\cdot g(x)$. A functional  $E'(x_*)$  we call a  unit 
bounded  root  functional  of polynomials  $f(x)$,  if 
  it   annuls $(f(x))^{\leq {\delta} _f+{\varepsilon} } _x$, 
where  ${\varepsilon} \geq 0$,  and  $E'(x_*)*1 = E'(y_*).\det 
\left\| \nabla f(x,y)\right\|= 1+f(x)\cdot g(x)$. 

{\it {\bf Theorem 2.}  Let ${\bf R} $   be a field. 
Let $x=(x_1,\ldots ,x_n)$   be  variables,  
let $f(x)=(f_1(x),\ldots ,f_n(x))$  be polynomials, 
${\delta} _f =\sum\limits^{ n} _{i=1} (\deg (f_i)-1)$. 
Let  ${\bf R} [x]/(f(x))_x$   be a finite-dimensional space 
over ${\bf R} $, in this case there exists a unit  root 
functional $E(x_*)$ of polynomials $f(x)$.\eject

 Let ${\varepsilon} \geq 0$, \ ${\cal A} (x_*)$  be the set of 
all functionals annulling $(f(x))^{\leq {\delta} _f+{\varepsilon} } _x$, 
${\cal L} (x_*)$ be the set of all functionals annulling 
$(f(x))_x$,  ${\cal U} (x)  =  {\bf R} [x^{\leq {\delta} _f+{\varepsilon} 
} ]$. Then:

 1) ${\cal A} (x_*)| _{{\cal U} (x)} $ with the extension operations 
for functionals 
is an associative and commutative  algebra 
over ${\bf R} $;

 2) there exists $d$ such that ${\cal A} (x_*)^d| _{{\cal U} 
(x)}   =  {\cal A} (x_*)^{d+1} | _{{\cal U} (x)} $,  and  for  
any such $d$ there holds  \\ \medskip

\ \ \ ${\cal L} (x_*)| _{{\cal U} (x)} 
 = {\cal A} (x_*)^d| _{{\cal U} (x)}  = ({\cal A} (x_*)| _{{\cal 
U} (x)} )^d.$} \medskip \\

 {\bf Proof  1.}   By virtue of  {\it  2}  of theorem  1  
 the extension operation for functionals 
induces the map \\ \medskip

\ \ \ $((f(x))^{\leq 
{\delta} _f+{\varepsilon} } _x)^\bot | _{{\bf R} [x^{\leq {\delta} 
_f+{\varepsilon} } ]} \times ((f(x))^{\leq {\delta} _f+{\varepsilon} } _x)^\bot 
| _{{\bf R} [x^{\leq {\delta} _f+{\varepsilon} } ]} 
 \rightarrow  ((f(x))^{\leq {\delta} _f+{\varepsilon} } _x)^\bot 
| _{{\bf R} [x^{\leq {\delta} _f+{\varepsilon} } ]}, $ 
\\ \medskip
\\
since for ${\delta} _1 = {\varepsilon} $, ${\delta} _2 = {\varepsilon} $, 
${\delta} = {\varepsilon} $, 
${\delta} '_1 = {\varepsilon} $, ${\delta} '_2 = {\varepsilon} $,  
${\delta} ' = {\varepsilon} $ 
there holds the condition {\it 2} of this theorem. Hence, ${\cal 
A} (x_*)| _{{\cal U} (x)}  = 
((f(x))^{\leq {\delta} _f+{\varepsilon} } _x)^
\bot | _{{\bf R} [x^{\leq {\delta} _f+{\varepsilon} }  ]} $  is  
an algebra with the extension operation 
for functionals.

 If $L_1(x_*), L_2(x_*) \in  {\cal A} (x_*)$, then they annul 
$(f(x))^{\leq {\delta} _f+{\varepsilon} } _x$. Then by virtue of 
theorem 1 in [7] 
$L_1(x_*)*L_2(x_*) = L_2(x_*)*L_1(x_*)$ in 
${\bf R} [x^{\leq {\delta} _f+{\varepsilon} +{\varepsilon} +1} ] \supseteq 
 {\bf R} [x^{\leq {\delta} _f+{\varepsilon} } ]$, and so, and 
in ${\bf R} [x^{\leq {\delta} _f+{\varepsilon} } ]$. 
This implies the commutativity of ${\cal A} (x_*)| _{{\cal U} (x)} $.

 If $L_1(x_*), L_2(x_*), L_3(x_*) \in  {\cal A} (x_*)$, then 
they annul $(f(x))^{\leq {\delta} _f+{\varepsilon} } _x$.  
Then by virtue of {\it 1} of theorem 2 in [7]    
$(L_1(x_*)*L_2(x_*))*L_3(x_*) = L_1(x_*)*(L_2(x_*)*L_3(x_*))$  
in ${\bf R} [x^{\leq {\delta} _f+{\varepsilon} +
{\varepsilon} +{\varepsilon} +2} ]  \supseteq 
  {\bf R} [x^{\leq {\delta} _f+{\varepsilon} } ]$,  
and so, and in 
${\bf R} [x^{\leq {\delta} _f+{\varepsilon} } ]$. This implies 
the associativity of ${\cal A} (x_*)| _{{\cal U} (x)} $.

 {\bf Proof 2.} In papers [1,3,4,5] 
there is 
the theorem 
about  existence of a unit root functional of polynomials $f(x)$ 
in the case, when ${\bf R} [x]/(f(x))_x$  be a finite-dimensional space 
over ${\bf R} $.

 An functionals in ${\cal A} (x_*)$ annul 
$(f(x))^{\leq {\delta} _f+{\varepsilon} } _x$, 
then by virtue of {\it 2} of theorem 3 in [6] 
any functional 
$L'(x_*) \in  {\cal A} (x_*)^p$  annuls 
$(f(x))^{\leq {\delta} _f+p\cdot {\varepsilon} +(p-1)} _x$, 
and so, annuls  
$(f(x))^{\leq {\delta} _f+p\cdot {\varepsilon} +(p-1)} _x\cap 
{\bf R} [x^{\leq {\delta} _f+{\varepsilon}}]$. \ 
By \ the \ finite \ dimensionality of 
${\bf R} [x^{\leq {\delta} _f+{\varepsilon} } ]$ 
\ over \ ${\bf R} $, \ there \ exists \ such  $p$, \ that \ 
$(f(x))^{\leq {\delta} _f+p\cdot 
{\varepsilon} +(p-1)} _x\cap {\bf R} [x^{\leq {\delta} _f+{\varepsilon} } ]  
=(f(x))_x\cap {\bf R} [x^{\leq {\delta} _f+{\varepsilon} } ]$, 
denote by $d=p$. Hence,  any  functional  
$L'(x_*)  \in {\cal A} (x_*)^d$ annuls 
$(f(x))_x\cap {\bf R} [x^{\leq {\delta} _f+{\varepsilon} } ]$. 
Then by virtue of {\it 4}  of theorem  6  in  [7] 
the functional $L(x_*) = L'(x_*)*E(x_*)$ 
annuls $(f(x))_x$ and $L(x_*) = L'(x_*)$ 
in ${\cal U} (x)  ={\bf R} [x^{\leq {\delta} _f+{\varepsilon} 
} ]$. Since $L'(x_*)$  is an arbitrary element $\in  {\cal 
A} (x_*)^d$, and $L(x_*) \in   {\cal L} (x_*)$, then ${\cal A} (x_*)^d| 
_{{\cal U} (x)}  \subseteq  {\cal L} (x_*)| _{{\cal U} (x)} $.

 Let $L(x_*)$ is an arbitrary element $\in  {\cal L} (x_*)$, 
then $L(x_*)$ annuls $(f(x))_x$. By virtue of {\it 2} of theorem 
6 in [7] $L(x_*)*E(x_*)  =  L(x_*)$.  Since  and  $E(x_*)  
\in {\cal L} (x_*)$, then 
\hbox{${\cal L} (x_*)*{\cal L} (x_*) = {\cal L} (x_*)$,} 
and so, ${\cal L} (x_*)^d = {\cal L} (x_*)$. 
There holds  ${\cal L} (x_*)  \subseteq {\cal A} (x_*)$, since 
any functional, annulling $(f(x))_x$, annuls $(f(x))^{\leq 
{\delta} _f+{\varepsilon} } _x$. Hence, ${\cal L} (x_*) 
= {\cal L} (x_*)^d \subseteq  {\cal A} (x_*)^d$, and so, 
${\cal L} (x_*)| _{{\cal U} (x)}  \subseteq  
{\cal A} (x_*)^d| _{{\cal U} (x)} $.

 From above it follows that ${\cal A} (x_*)^d| _{{\cal 
U} (x)}  = {\cal L} (x_*)| _{{\cal U} (x)} $.

 Since by virtue of {\it 1} of the theorem 
the extension map $*$ for functionals induces the map 
\hbox{${\cal A} (x_*)| _{{\cal U} (x)} \times {\cal A} (x_*)| _{{\cal U} (x)}$} 
$\rightarrow 
{\cal A} (x_*)| _{{\cal U} (x)} $,
then $({\cal A} (x_*)| _{{\cal U} (x)} )^d={\cal A} (x_*)^d| _{{\cal U} (x)} $. \eject

 {\bf Algorithm.}  (Finding a basis  of all root functionals 
and  a basis of the ideal of polynomials, and also the unit root functional.) 
Let ${\bf R} $   be a field,  let $x=(x_1,\ldots ,x_n)$, $y\simeq 
x$  be  variables,  $f(x)=(f_1(x),\ldots ,f_n(x))$   be  polynomials. 
 Let ${\bf R} [x]/(f(x))_x$  be a  finite-dimensional   space 
  over   ${\bf R} $.   Denote by   ${\delta} _f   =\sum\limits^{ 
n} _{i=1} (\deg (f_i)-1)$,  ${\cal L} (x_*)  =  ((f(x))_x)^\bot 
$,  ${\cal U} (x)={\bf R} [x^{\leq {\delta} _f} ]$.  Here  and 
 below  by space we shall mean a linear space over ${\bf R} $. 

 The algorithm finding a basis of the space  
${\cal L} (x_*)| _{{\cal U} (x)}   =  
((f(x))_x)^\bot | _{{\bf R} [x^{\leq {\delta} _f} ]} $ 
of restrictions of all  root  functionals  on 
${\bf R} [x^{\leq {\delta} _f} ]$  
and  a basis of the space 
$(f(x))_x\cap {\bf R} [x^{\leq {\delta} _f} ]$, 
and also the restriction of the unit root functional on 
${\bf R} [x^{\leq {\delta} _f} ]$ consists of the following  steps:

 1. Construct by Gauss elimination method a basis of the space 
$(f(x))^{\leq {\delta} _f} _x$.

 2. From Gauss basis of the space $(f(x))^{\leq {\delta} _f} _x$  
construct  Gauss  basis of the space of functionals 
\hbox{defined} on ${\bf R} [x^{\leq {\delta} _f} ]$ and annulling 
 $(f(x))^{\leq {\delta} _f} _x$,  this space coincide 
with ${\cal A} (x_*)| _{{\cal U} (x)}  = ((f(x))^{\leq {\delta} 
_f} _x)^\bot | _{{\bf R} [x^{\leq {\delta} _f} ]} $. Let  obtained 
basis be  $L_1(x_*),\ldots ,L_d(x_*)$.

 3. Compute the restriction of operators 
$\left[ L_1(x_*)\right] ,\ldots ,\left[ L_d(x_*)\right] $ 
on ${\bf R} [x^{\leq {\delta} _f} ]$.

 4.  Compute  the restriction of functionals   $(L_1(x_*))^d,\ldots 
,(L_d(x_*))^d$   on ${\bf R} [x^{\leq {\delta} _f} ]$ by 
\smallskip \\

\ \ \ $ \forall p=1,d: \forall {\delta} =2,d: (L_p(x_*))^{\delta} 
 = (L_p(x_*))^{{\delta} -1} .\left[ L_p(x_*)\right]  \hbox{ in } {\bf R} [x^{\leq {\delta} _f} ].$
\smallskip \\

5. Compute  the following  generators:  \\ \smallskip

\ \ \ $\{(L_p(x_*))^d*L_q(x_*)|p=1,d \  \& \  q=1,d\}$ 
\\ \smallskip
\\
of the space 
${\cal L} (x_*)| _{{\cal U} (x)} =({\cal A} (x_*)| _{{\cal U} (x)} )^{d+1}.$

 6. By Gauss elimination method construct a basis of the space 
 ${\cal L} (x_*)| _{{\cal U} (x)} $ from its system of generators.

 7. From Gauss basis of the space  ${\cal L} (x_*)| _{{\cal 
U} (x)} $  construct  Gauss  basis of the  space of polynomials  
$\in   {\cal U} (x)  =  {\bf R} [x^{\leq {\delta} _f} ]$  annulled by 
  ${\cal L} (x_*)| _{{\cal U} (x)} $.  This space of polynomials 
coincide with 
\hbox{$(f(x))_x\cap {\bf R} [x^{\leq {\delta} _f} ]$.}

 8. Let $h_1(x),\ldots ,h_{d\,'} (x)$  be a basis of the space 
 $(f(x))_x\cap {\bf R} [x^{\leq {\delta} _f} ]$,  
let $l_1(x_*),\ldots ,l_{d\,''} (x_*)$ 
be a basis of the space 
${\cal L} (x_*)| _{{\cal U} (x)} $. 
From $\{l_p(y_*).\det \left\| \nabla f(x,y)\right\| | p=1,d\,''\}$  
and  $\{h_q(x)| q=1,d\,'\}$  
by Gauss elimination method find the decomposition   
\\
$$\sum\limits^{ d\,''} _{p=1} a_p\cdot (l_p(y_*).
\det \left\| \nabla f(x,y)\right\| )  +   
\sum\limits^{ d\,'} _{q=1} b_q\cdot h_q(x)   =   1,
\hphantom{ccccccccccccccccccccccccccccccccc} $$ 
\\
$E'(x_*)=\sum\limits^{ d\,''} _{p=1} a_p\cdot l_p(x_*)$ is  the restriction of  the unit 
 root  functional of polynomials $f(x)$ 
on ${\bf R} [x^{\leq {\delta} _f} ]$,  
since  $E'(y_*).\det \left\| \nabla f(x,y)\right\| - 1 \in (f(x))_x$. \eject

 {\bf Proof of the algorithm.}

 5. The dimension of the space ${\cal A} (x_*)| _{{\cal U} (x)} 
$ is equal to $d$. Therefore the chain \medskip

\ \ \ $({\cal A} (x_*)| _{{\cal U} (x)} 
)^1 \supseteq  \ldots  \supseteq  ({\cal A} (x_*)| _{{\cal U} 
(x)} )^{\delta}  \supseteq  ({\cal A} (x_*)| _{{\cal U} (x)} 
)^{{\delta} +1}  \supseteq  \ldots $ 
\medskip 
\\
is stabilized for 
some 
${\delta} \leq d+1$, i. e. 
$({\cal A} (x_*)| _{{\cal U} (x)} )^{{\delta} '} = 
({\cal A} (x_*)| _{{\cal U} (x)} )^{{\delta} '+1} $ for 
 any ${\delta} '\geq {\delta} $. 
\smallskip Then 

\ \ \ $({\cal A} (x_*)| 
_{{\cal U} (x)} )^{d+1}  = ({\cal A} (x_*)| _{{\cal U} (x)} )^{d+2} 
 = \ldots  = ({\cal A} (x_*)| _{{\cal U} (x)} )^{d+{\delta} '} 
 = \ldots \ . $ 
\medskip \\
Any element in $({\cal A} (x_*)| _{{\cal U} 
(x)} )^{d+1}  = ({\cal A} (x_*)| _{{\cal U} (x)} )^{d\cdot 
(d-1)+2} $   
generated by 
elements of the form \ $(L_1(x_*))^{{\alpha} 
_1} *\ldots *(L_d(x_*))^{{\alpha} _d} $, where 
${\alpha} _1+\ldots +{\alpha} _d = d\cdot (d-1)+2$.  
Then there exists  such  $p$,   that  ${\alpha} _p\geq d$,  
since otherwise  $\forall p:  {\alpha} _p\leq d-1$,  
and, hence, 
$d\cdot (d-1) \geq  {\alpha} _1+\ldots +{\alpha} _d = d\cdot (d-1)+2$, 
that is impossible.  In this case 

\ \ \ $(L_1(x_*))^{{\alpha} _1} *\ldots *(L_d(x_*))^{{\alpha} _d}  
= (L_p(x_*))^d*\left( (L_p(x_*))^{{\alpha} _p-d} *(\prod\limits_{q:\not= p} 
(L_q(x_*))^{{\alpha} _q} )\right) =$
  
\ \ \ $\qquad = (L_p(x_*))^d*L(x_*). $ \medskip
\\
Here  $L(x_*)  \in   {\cal A} (x_*)| _{{\cal 
U} (x)} $,  since  $({\alpha} _1+\ldots +{\alpha} _d)  - 
 d  =d\cdot (d-1)+2 - d = d\cdot (d-2) +2 = (d-1)^2+1 \geq 
 1$. Then $L(x_*)$ is expressed via $L_1(x_*),\ldots 
,L_d(x_*)$   linearly  over  ${\bf R} $.   Hence, the  space 
$({\cal A} (x_*)| _{{\cal U} (x)} )^{d+1} $ is generated by generators 
\medskip
 
\ \ \ $\{(L_p(x_*))^d*L_q(x_*)| p=1,d \ \& \ q=1,d\}.$ 
\medskip
\\
That ${\cal L} (x_*)| _{{\cal U} (x)}  = ({\cal A} (x_*)| _{{\cal U} (x)} 
)^{d+1} $ is stated in {\it 2} of theorem 2.

 7. Any functional in 
${\cal L} '(x_*) = {\cal L} (x_*)| _{{\cal U} (x)}  = 
((f(x))_x)^\bot | _{{\bf R} [x^{\leq {\delta} _f} ]} $ 
annuls the space 
${\cal M} (x) =  (f(x))_x\cap {\bf R} [x^{\leq {\delta} _f} ]\subseteq 
{\cal U} (x)$,  
hence, 
${\cal L} '(x_*)  \subseteq   {\cal M} (x)^\bot $. Here  we consider 
the annulet of the space ${\cal M} (x)$   
as a subspace of the space 
${\cal U} (x)$,  and the annulet of the space  
${\cal L} '(x_*)$ as a subspace of the space 
${\cal U} (x)_*$. 
Let a functional $l(x_*)$, 
determined on  ${\cal U} (x) =  {\bf R} [x^{\leq {\delta} _f} ]$, 
annuls 
${\cal M} (x) = (f(x))_x\cap {\bf R} [x^{\leq {\delta} _f} ]$. 
By virtue of statement 2 functional  
$L'(x_*)  =l(x_*).{\cal P} ^{\leq {\delta} _f} _x$ 
is determined on ${\bf R} [x]$,  and  
$L' (x_*)| _{{\bf R} [x^{\leq {\delta} _f} ]} = l(x_*)$,  
hence, 
$L'(x_*)$ annuls 
${\cal M} (x) = (f(x))_x\cap {\bf R} [x^{\leq {\delta} _f} ]$. 
Then by virtue of {\it 4} of theorem  6 in [7] 
there exists $L(x_*)$, annulling 
$(f(x))_x$, such that 
$L(x_*)| _{{\cal U} (x)}  = L(x_*)| _{{\bf R} [x^{\leq {\delta} _f} ]}=
L'(x_*)| _{{\bf R} [x^{\leq {\delta} _f} ]}= l(x_*)$. Since $L(x_*) 
\in {\cal L} (x_*)$, 
then 
$l(x_*) \in {\cal L} (x_*)| _{{\cal U} (x)}  
 ={\cal L} '(x_*)$. 
Hence, 
${\cal M} (x)^\bot \subseteq {\cal L} '(x_*)$. 
Thus  ${\cal L} '(x_*) = {\cal M} (x)^\bot $. 
Then  by virtue of the finite dimensionality of the space 
${\cal U} (x) = {\bf R} [x^{\leq {\delta} _f} ]$ 
there holds ${\cal M} (x)  =  {\cal L} '(x_*)^\bot $. Here we identify 
$({\cal U} (x)_*)_*$ with ${\cal U} (x)$. 

 {\bf Estimation of the complexity of the algorithm.}  Let $D$  be a dimension  
of the space  ${\bf R} [x^{\leq {\delta} _f} ]$, then  $D  
=  C^n_{{\delta} _f+n}   =  C^n_{d_1+\ldots +d_n} $,  where  $d_i=\deg 
(f_i)$.  Let us estimate  the complexity of steps of the algorithm.

 1. The number of polynomials in system of polynomials $\{f_i(x)\cdot x^{{\alpha} 
(i)} \in  {\bf R} [x^{\leq {\delta} _f} ]|   i=1,n\}$  not exceed 
$n\cdot D$. Construction of a basis  from  this  system  of polynomials  
by Gauss elimination method  requires 
\hbox{$\leq  (n\cdot D)\cdot O(D^2) = n\cdot O(D^3)$} 
operations.

 2. The step  requires $\leq  O(D^2)$ operations. \eject

 3. Computation of all minors of the matrix \medskip 

\ \ \ $\left\| \begin{matrix}\nabla f(x,y) \cr
 f(x)\end{matrix} \right\| $ 
\medskip
\\
of order $n$ without divisions  
 requires  $\leq   ({\delta} _f\cdot n^2+n^4)\cdot O(D^3)$  operations. 
Within this 
it computing and $\det \left\| \nabla f(x,y)\right\| 
$. Computation of the operator \medskip

\ \ \ $\left[ L_p(x_*)\right]  = L_p(y_*).\det 
\left\| \begin{matrix}\nabla f(x,y) & \nabla _x(x,y) \cr
 f(x) &  {\bf 1} _x(x)\end{matrix} \right\| $
\medskip
\\
on ${\bf R} 
[x^{\leq {\delta} _f} ]$  requires $\leq  O(D^3)$ operations. 
Computation of  such  operators  for all $p=1,d$  requires 
$\leq  d\cdot O(D^3) \leq  D\cdot O(D^3) = O(D^4)$ operations.

 4. Computation of $(L_p(x_*))^{\delta}  = (L_p(x_*))^{{\delta} 
-1} .\left[ L_p(x_*)\right] $    requires  $\leq   O(D^2)$ 
operations, and for all ${\delta} =2,d$ and $p=1,d$  requires 
 $\leq   d^2\cdot O(D^2)  \leq   D^2\cdot O(D^2)  =O(D^4)$ 
operations.

 5. Computation of $(L_p(x_*))^d*L_q(x_*) = (L_p(x_*))^d.\left[ 
L_q(x_*)\right] $  requires $\leq O(D^2)$ operations. Since 
this computation necessary to perform for all $p=1,d$ and of all 
$q=1,d$, then in all performed  $\leq  d^2\cdot O(D^2) \leq 
 D^2\cdot O(D^2) = O(D^4)$ operations.

 6. Computation of a basis of the space ${\cal L} (x_*)| _{{\cal U} (x)} $  
from  $d^2$  generators by Gauss elimination method  
requires $\leq  d^2\cdot O(D^2) \leq  D^2\cdot O(D^2) = O(D^4)$ operations.

 7.  The step   requires $\leq  O(D^2)$ operations. 

 8. Computation of $l_p(y_*).\det \left\| \nabla f(x,y)\right\| 
$ for  single  $p$    requires  $\leq O(D^2)$ operations, 
then computation for all $p=1,d''$  requires 
$\leq  d''\cdot O(D^2) \leq  D\cdot O(D^2)= O(D^3)$ 
operations. Decomposition of $1$ by Gauss elimination method requires 
$\leq  O(D^3)$ operations, and computation of 
$E'(x_*) = \sum\limits^{ d''} _{p=1} a_p\cdot l_p(x_*)$  
requires $\leq  d''\cdot O(D) \leq  D\cdot O(D)  \leq   O(D^2)$ operations.

 If to regard $n$ as constant, then  the summarized  number of operations 
 performed  in the algorithm is $\leq O(D^4)$. \medskip  {\it

 {\bf Theorem 3.}  Let $x=(x_1,\ldots ,x_n)$   be  variables, 
 $f(x)=(f_1(x),\ldots ,f_n(x))$ be polynomials, 
${\delta} _f = \sum\limits^{ n} _{i=1} (\deg (f_i)-1)$. 
Let ${\bf R} [x]/(f(x))_x$  
be a finite generated as module over ${\bf R} $, then  
\medskip 

\ \ \ $(f(x))_x\cap {\bf R} [x^{\leq {\delta} _f+{\delta} } ] = 
((f(x))_x\cap {\bf R} [x^{\leq {\delta} _f} ])\cdot 
{\bf R} [x^{\leq {\delta} } ] + 
(f(x))^{\leq {\delta} _f+{\delta} } _x,$ \medskip

\ \ \ $(f(x))_x\cap {\bf R} [x^{\leq {\delta} _f+{\delta} +{\varepsilon} +1} ] = 
((f(x))_x\cap {\bf R} [x^{\leq {\delta} _f+{\varepsilon} +1} ])\cdot 
{\bf R} [x^{\leq {\delta} } ],$  \medskip
\\
where  ${\delta} 
\geq 0$  and ${\varepsilon} \geq 0$. } \smallskip

 Proof of theorem 3 will be given in the subsequent papers. 
\medskip

{\footnotesize

\begin{enumerate}

\item {\it Seifullin, T. R.} 
Root functionals and root polynomials of a system of polynomials. (Russian)
Dopov. Nats. Akad. Nauk Ukra\"\i ni  -- 1995, -- no. 5, 5--8.
\item {\it Seifullin, T. R.} Root functionals and root relations 
of a system of polynomials. (Russian) 
Dopov. Nats. Akad. Nauk Ukra\"\i ni  -- 1995, -- no 6, 7--10.
\item  {\it Seifullin, T. R.}  Homology of the Koszul complex of a 
system of polynomial equations. (Russian)
Dopov. Nats. Akad. Nauk Ukr. Mat. Prirodozn. Tekh. Nauki 1997, no. 9, 43--49. 
\item  {\it Seifullin, T. R.}  Koszul complexes of systems of 
polynomials connected by linear dependence. (Russian) 
Some problems in 
contemporary mathematics (Russian), 326--349, Pr. Inst. Mat. Nats. Akad. Nauk 
Ukr. Mat. Zastos., 25, Natsional. Akad. Nauk Ukra\"\i ni, Inst. Mat., Kiev, 
1998. \eject
\item {\it Seifullin, T. R.} Koszul complexes of embedded systems of 
polynomials and duality. (Russian) 
Dopov. Nats. Akad. Nauk Ukr. Mat. Prirodozn. 
Tekh. Nauki 2000, no. 6, 26--34.  
\item {\it Seifullin, T. R.} 
Extension of bounded root functionals 
of a system of polynomial equations. Dopov. Nats. Akad. Nauk Ukr. Mat. 
Prirodozn. Tekh. Nauki 2002, no. 7, 35--42. 
\href{http://arxiv.org/abs/0804.2420}
{{\tt  arXiv:0804.2420}}.
\item {\it Seifullin, T. R.} 
Continuation of root functionals of 
a system of polynomial equations and the reduction of polynomials modulo its 
ideal. (Russian) Dopov. Nats. Akad. Nauk Ukr. Mat. Prirodozn. Tekh. Nauki 2003, 
no. 7, 19--27.  
\href{http://arXiv.org/abs/0805.4027}
{{\tt  arXiv:0805.4027}} (English).
\item {\it  Buchberger B.} Gr\"obner: An algorithmic method in polynomial
ideal theory
//Multidimensional Systems Theory. / Ed. N. K. Bose, --
Dordrecht: D. Reidel, 1985. -- Chapter 6.
\item {\it Caniglia L., Galligo A., Heintz J.} Some new effictivity bounds
in computational geometry
// Proc. 6th Int. Conf. on Appied Algebra and Error--correcting codes.
/ LNCS 357, Springer--Verlag, Berlin. -- 1989. -- pp. 131--152.
\item {\it Brownawell D.} Bounds for the degrees in the Nullstellensatz 
// Ann. Math. 2nd series. -- 1987. -- No 126. -- pp. 577--591.
\item {\it Canny J.} Generalized characteristic polynomials 
//J. Symbolic Computation. -- 1990. -- No 9. -- pp. 241--250.
\\
\end{enumerate}
\small{\noindent
{\it V. M. Glushkov Institute of Cybernetics of the NAS of Ukraine, Kiev
\hfill Received 26.06.2002\medskip\\
E-mail: \ {\tt  timur\_sf@mail.ru}
}

\end{document}